\documentstyle[11pt,amstex,righttag]{amsart}

\theoremstyle{definition}

\theoremstyle{remark}

\title[On the average rate of return ...]{ On the average
rate of return in a continuous time stochastic
model${}^{\bf 1}$}
\date{December 31, 2000}
\thanks{${}^{1}$ Supported by the KBN Grant  No. 1H02B 018 14.}
\author[L.Gajek and M.Kaluszka]{Leslaw Gajek, Marek Kaluszka}
\address{Institute of Mathematics\\ Technical University of Lodz\\ ul.
Zwirki 36\\ 90-924 Lodz\\ Poland} 

\begin{document}
\maketitle

\begin{abstract}
In a discrete time stochastic model of a pension investment funds market
Gajek and Kaluszka(2000a) have provided a definition of the
average rate of return which satisfies a set of economic
correctnes postulates. In this paper the average rate of return is
defined for a continuous time stochastic model of the market. The prices of
assets are modeled by the multidimensional geometrical Brownian motion.
A martingale property of the average rate of return is proven.

\end{abstract}
\medskip
\medskip

\section{Introduction}

Consider a group of $n$ investment or pension funds.
In Gajek and Kaluszka (2000a)
the following definition of the average return of the
group of funds during a given time period $[s,t]$ is
proposed
$$
\bar{r}(s,t)=\exp\left(\sum^{n}_{i=1} \int^{t}_{s}
{ k_{i}(u)w_{i}(u)\over \sum^{n}_{j=1}k_{j}(u)w_{j}(u)}
\delta _{i}(u)du\right)-1,
\eqno{(1.1)}
$$
where $k_{i}(t)$ and $w_{i}(t)$ denote a number of all units possessed by the
members of the $i$-th fund and a value of the $i$-th fund unit at time $t$,
respectively, and
$\delta _{i}(t)$ is
the instantaneous rate of interest for the accounting unit
in the $i$-th fund defined as
$$
\delta _{i}(t)={d\over dt}\log w_{i}(t)
\eqno{(1.2)}
$$
Here and subsequently, the logarithm is to the base $e$. This
definition possesses some desirable properties contrary to
other definitions (see Gajek and Kaluszka (2000a))

In Section 2 we introduce a definition of the average rate of return
in a continuous time stochastic model based on a
discrete time model presented in Gajek and Kaluszka(2000b). We assume
that prices of assets are modelled as correlated geometric Brownian motions.
This allows us to derive a proper definition of the average rate of return.
The definition is given by
$$ \bar{r}(s,t)=\exp \left(\sum^{n}_{i=1}
\int^{t}_{s} A_i^*(u) {dw_{i}(u) \over w_{i}(u)} +\int^{t}_{s}d
\log A(u) - \int^{t}_{s}{dA(u)\over A(u)}+\right.
$$
$$
\left.+\sum _{i=1}^n\int _s^t {A^* _i }^2(u){dk_i(u)\over k_i(u)}-
\sum _{i=1}^n\int _s^t {A^*_i }^2(u)d\log k_i(u)
\right)-1. \eqno{(1.3)}
$$
Here
and subsequently,
$$ A_{i}(t)=k_i(t)w_i(t), \quad A(t) =
\sum^{n}_{i=1}A_{i}(t), \quad A^*_i (u)=\frac{A_i(u)}{A(u)},\eqno{(1.4)}$$
and the intergals in  (1.3) are the It\^{o} integrals.

The suitability of the
formula is investigated.

\medskip
\section{A continuous-time stochastic model for funds dynamics}

\noindent {\bf 2.1. The model.}
We consider the following state-variables:

$u_{ij}(t) =$ number of units of the $j$-th asset possessed by the $i$-th
fund at time $t$, $i=1,\ldots  ,n$, $j=1,\ldots , N$.

$w_{i}(t) =$ value of a participation unit of the $i$-th fund at time $t$,

$k_{i}(t) =$  number of units of the $i$-th fund at time $t$,

$D_{i}(t)=$  inflow of contribution income minus an outflow
of benefit payments of the $i$-th fund over the time period $[0,t]$,

$D(t)=\sum^{n}_{i=1}D_{i}(t)$,

$K_{ij}(t)$ = number of units which have flowed from the $i$-th fund to $j$-th one between
time $0$ and time $T$, where $j\neq i$.

\medskip
Let $(\Omega ,{\cal F},{\Bbb P})$ be a complete probability space. Let
${\Bbb F}=\{{\cal F}_{t}\}$ be a filtration,
i.e. each ${\cal F}_{t}$ is an $\sigma $-algebra of subsets of
$\Omega $ with ${\cal F}_{s}\subseteq {\cal F}_{t}\subseteq {\cal F}$
for every $s<t$. Without loss of generality, we assume
that ${\cal F}_{0}=\{\Omega ,\emptyset \}$.
The investors' planning horizont is $T$, a fixed
positive number.

Assume we have $N$ assets and their prices $c_i$, $i=1,\ldots ,N$,  are
governed by the following stochastic differential equations: $$
dc_{i}(t)=c_{i}(t)\mu _{i}(t)dt+c_{i}(t)\sum^{N}_{j=1} \sigma
_{ij}(t)dB_{j}(t)\quad \hbox{ for }\quad i=1,\ldots  , N, \eqno{(2.1)} $$
where $(B_{1}(t), \ldots  , B_{N}(t))$ is a standard $N$-dimensional
Brownian motion (under the real-world probability ${\Bbb P})$, and $\sigma
_{ij}(t)$ and $\mu _{i}(t)$ are progressively measurable processes on
$[0,T]$ with bounded variation, called volatility processes and mean return
rate processes, respectively. We assume the matrix $(\sigma
_{ij})_{i,j=1}^N$ is non-singular and assume that $\sum _{i=1}^N\sum
_{j=1}^N\int^{T}_{0}\sigma _{ij}^2(t)dt<\infty $. Moreover, let us assume
that $\sum _{i=1}^N \int _0^T |\mu _i(t)|dt<\infty$. Here and subsequently,
the symbol $X=Y$ (resp. $X<Y$) means that the random variables $X$, $Y$ are
defined on $(\Omega ,{\cal F},{\Bbb P})$ and ${\Bbb P}(X=Y)=1$ (resp.
${\Bbb P}(X<Y)=1$). The model (2.1) of prices of assets is commonly used in
mathematical finance (see e.g. Musiela and Rutkowski(1997), Karatzas and
Shreve(1998) or Shiryaev(1999)). By It\^{o}'s formula $$
c_{i}(t)=c_{i}(0)\exp \left(\int^{t}_{0} \left[\mu _{i}(s)-\frac{1}{2}\sum
_{k=1}^N \sigma _{ik}^2(s)\right]ds+ \sum^{N}_{k=1} \int^{t}_{0} \sigma
_{ik}(s)dB_{k}(s)\right), $$ for $i=1,\ldots  , N$.

The process $D_i$ will be modelled as follows
$$
dD_i(t)=\alpha _i(t)dt+\beta _i(t)dB_{i+N}(t),\quad i=1,\ldots , N,\eqno{(2.2)}
$$
where $\alpha _i$, $\beta _i$ are progressively measurable processess,
and $\{B_{i}(t):\; i=1,2,\ldots , 2N\}$  are independent standard Brownian motions.
The class (2.2) contains geometric Brownian motion which is commonly used
as a model of inflow-outflow process (cf. Koo(1998)).

We assume that all investments are infinitely divisible. There are
no transaction costs or taxes, the assets pay no dividends and
member does not pay for allocation of his/her wealth. The split of
units is not allowed.  In addition,
suppose that $A_i(t)>0$ between times $0$ and $T$ for each $i$.
Assuming strict positivity
we avoids various technical complications.

The dynamics of a group of funds is described by the following
stochastic differential equations:
$$
w_{i}(t)k_{i}(t)=u_{i1}(t)c_{1}(t)+\ldots  +u_{iN}(t)c_{N}(t),
\eqno{(2.3)}
$$
$$
k_{i}(t)dw_{i}(t)=u_{i1}(t)dc_{1}(t)+\ldots  +u_{iN}(t)dc_{N}(t),
\eqno{(2.4)}
$$
$$
w_{i}(t)dk_{i}(t)=-w_{i}(t)\sum^{}_{j\neq i}
dK_{ij}(t)+\sum^{}_{j\neq i}w_{j}(t)dK_{ji}(t)+dD_{i}(t),
\eqno{(2.5)}
$$
$$
w_{i}(t)dk_{i}(t)=c_{1}(t)du_{i1}(t)+\ldots
+c_{N}(t)du_{iN}(t),
\eqno{(2.6)}
$$
where $i=1,2,\ldots  ,n$, $t\in [0,T]$, and processes $c_i$, $D_i$ are
defined by (2.1) and (2.2), respectively. The equations (2.3)-(2.6) are
counterparts of equations (2.8) and (2.10)-(2.12) of Gajek and Kaluszka(2000b).
The functions
$u_{ij}$, $K_{ij}$ play a role of control variables.
We also assume that $K_{ij}$ is continuous stochastic process adapted to
${\Bbb F}$
with bounded variation on any compact interval.

After adding equations (2.5) we get
$$
\sum^{n}_{i=1}w_{i}(t)dk_{i}(t)=d D(t).
\eqno{(2.7)}
$$

{\bf Remark.} In Bacinello(2000)
(see also Chamorro and de Villarreal(2000)),
the value $w_i(t)$ of participation unit of the $i$-th fund is
modelled as geometric Brownian motion.\; In our setting, each random variable
$w_i(t)$ is a sum of lognormaly distributed random variables. It seems
that it is not possible to derive a proper
definition of the average rate of return under the assumption of
Bacinello(2000). \medskip

Recall It\^{o}'s formula for
$f(\xi )$ with $\xi =\xi (t)$ such that
$d\xi =\alpha dt+\sum \beta _{j}dB_{j}$:
$$
df(\xi )=f'(\xi )\sum \beta _{j}dB_{j}+
\left(f'{}(\xi )\alpha +{1\over 2}f'{}'{}(\xi )\sum \beta _{j}^2\right)dt
$$
$$
=f'(\xi )d\xi  +{1\over 2}f'{}'{}(\xi )\left(\sum \beta _{j}^2\right)dt,\eqno{(2.8)}
$$
where $f\in C^{2}(\Bbb R)$, and $\alpha $, $\beta _i$ are progressively measurable
processess (see also the  local It\^{o}'s formula, e.g. Kallenberg(1997)).

Moreover, if $\xi _i= \alpha dt+\sum _j\beta _{ij}dB_{j}$ for $i=1,2$, then
$$
d(\xi _1\xi _2)=\xi _1d\xi _2 + \xi _2d\xi _1+
(\sum \beta _{1j}\beta _{2j})dt.\eqno{(2.9)}
$$

\medskip
\noindent {\bf 2.2. Definition of the average return.} Our definition of the
average rate of return
of a group of funds at a time
interval $[s,t]$, say $\bar{r}(s,t)$,
in the stochastic continuous-time model, is as follows
$$ \bar{r}(s,t)=\exp \left(\sum^{n}_{i=1}
\int^{t}_{s} A_i^*(u) {dw_{i}(u) \over w_{i}(u)} +\int^{t}_{s}d
\log A(u) - \int^{t}_{s}{dA(u)\over A(u)}+\right.
$$
$$
\left.+\sum _{i=1}^n\int _s^t {A^*_i } ^2(u){dk_i(u)\over k_i(u)}-
\sum _{i=1}^n\int _s^t {A^*_i } ^2(u)d\log k_i(u)
\right)-1. \eqno{(2.10)}$$
Recall that $A_{i}(t)$ denotes the wealth of the $i$-th fund at time $t$,
$A(t)$
means the global wealth of funds at time $t$, and $A^{*}_i(t)$ denotes  the
percentage of a relative value of assets of the $i$-th fund.
The formula (2.10) is a generalization of (1.1) to the case of stochastic
prices and random inflow-outflow process
modelled by (2.1) and (2.2), respectively.
In fact, if $\beta _i\equiv 0$ , i.e. $dD_i(t)=\alpha _i(t)dt$, then
$$
\bar{r}(s,t)=\exp \left(\sum^{n}_{i=1}
\int^{t}_{s} A^{*}_{i}(u){dw_{i}(u)
\over w_{i}(u)} +\int^{t}_{s}d \log A(u)
- \int^{t}_{s}{dA(u)\over A(u)}\right)-1.
\eqno{(2.11)}
$$
Moreover, if $dc_i(t)=\mu _i(t)dt$ for each $i$, then
$$
\bar{r}(s,t)=\exp \left(\sum^{n}_{i=1}
\int^{t}_{s} A^{*}_{i}(u){dw_{i}(u)
\over w_{i}(u)}\right)-1.
$$
Let us now provide an heuristic argument for definition (2.11) based on the
definition of the average rate of return $\bar{r}_A$ in a discrete time model:
$$
\bar{r}_{A}(s,t) =\prod^{t-1}_{u=s}
\left(1+\sum^{n}_{i=1} A^{*}_{i}(u)r_{i}(u,u+1)\right) -1.
$$
proposed by Gajek and Kaluszka(2000b).  By
formula (2.15) of Gajek et al. (2000b)
$$
\bar{r}_{A}(s,t)=\prod _{u=s}^{t-1} {A(u+1)-d(u+1)\over A(u)} -1 $$
$$={A(t)\over A(s)}\exp\left(
\sum^{t-1}_{u=s} \log (1- {D(u+1)-D(u)\over A(u)})\right)-1 $$ $$ \approx
{A(t)\over A(s)}\exp\left( - \sum^{t-1}_{u=s}{D(u+1)-D(u)\over
A(u)}\right)-1 $$
$$ \approx  {A(t)\over
A(s)}\exp\left(-\int^{t}_{s}{dD(u)\over A(u)}\right)-1,\eqno{(2.12)}
$$ since
$\ln(1-x)\approx  -x$ if $x$ is close to 0. The latter approximation in (2.12) follows from
the definition of It\^{o}'s integral.
Combining (2.3) and (2.6) we get
$$ dD(u)= \sum^{n}_{i=1}w_{i}(u)dk_{i}(u)= dA(u)-
\sum^{n}_{i=1}k_{i}(u)dw_{i}(u). \eqno{(2.13)} $$
From (2.12) and (2.13) it follows
$$
 \bar{r}_{A}(s,t)\approx  {A(t)\over A(s)} \exp\left(\sum^{n}_{i=1}
\int^{t}_{s}k_{i}(u) {dw_{i}(u)\over A(u)} - \int^{t}_{s}{dA(u)\over
A(u)}\right)-1
$$
$$
={A(t)\over A(s)} \exp\left(\sum^{n}_{i=1}
\int^{t}_{s}A_{i}^*(u) {dw_{i}(u)\over w_i(u)} - \int^{t}_{s}{dA(u)\over
A(u)}\right)-1.
$$
Since  $A(t)$ is a strictly positive process,
$$
\bar{r}_A(s,t)=
\exp \left(\sum^{n}_{i=1} \int^{t}_{s} A^{*}_{i}(u){dw_{i}(u)
\over w_{i}(u)}+\int _s^td\log A(u) - \int^{t}_{s}{dA(u)\over A(u)}\right)-1
\eqno{(2.14)}
$$
(by the local It\^{o}'s formula), we obtain
$$
\bar{r}_A(s,t)\approx \bar{r}(s,t)
$$
as claimed.
\medskip

\noindent {\bf 2.3. Properties of the average rate of return (2.10)}

\medskip
{\bf Theorem 1.} Suppose $\sigma _{ij}(\cdot )$ are bounded in $t$ and $\omega $.
If $\{c_i(t), t\ge 0\}$ is an
${\Bbb F}$-martingale for every $i$, then $\{\bar{r}(0,t), t\ge 0\}$ is also
an ${\Bbb F}$-martingale.

{\it Proof.} To simplify notation
we cancel the dependence of time of considered processes.
Recall that $K_{ij}$ is adapted and continuous processes with
bounded variation on any compact interval for each $i,j$.
First observe that by (2.1), (2.2), (2.4) and (2.5) we have
$$
dw_i=\frac{1}{k_i}\left(\sum _{j=1}^N c_j\mu _j u_{ij}
\right)dt+
\frac{1}{k_i}\sum _{l=1}^N
\left(\sum _{j=1}^Nc_j\sigma _{jl} u_{ij}\right)dB_l,\eqno{(2.15)}
$$
$$
dk_{i}=\left(-w_i\sum^{}_{j\neq i}
k_{ij}+\sum^{}_{j\neq i}w_{j}k_{ji}+\alpha _i\right){dt\over w_i}
+{\beta _i\over w_i}dB_{i+N},\quad i=1,2,\ldots , N,
\eqno{(2.16)}
$$
where $dK_{ij}=k_{ij}dt$. By (2.3), (2.7), (2.15), (2.16) and It\^{o}'s
formula (see (2.9)) we get
$$
dA=\sum^{n}_{i=1}w_{i}dk_{i}+
\sum^{n}_{i=1}k_{i}dw_{i} + (0\cdot B_1+\ldots +0\cdot B_N+B_{N+1}\cdot 0+\ldots +B_{2N}\cdot 0)dt
$$
$$=dD+
\sum^{n}_{i=1}k_{i}dw_{i}=
(\sum _{i=1}^N\alpha _i +\sum^{N}_{j=1}c_{j}\mu _{j}\sum^{n}_{i=1}
u_{ij})dt+$$
$$+\sum^{N}_{l=1}(\sum^{N}_{j=1}c_{j}\sigma _{jl}
\sum^{n}_{i=1}u_{ij})dB_{l}+\sum _{l=1}^{N}\beta _ldB_{l+N}.\eqno{(2.17)}
$$
By local It\^{o}'s formula
we obtain
$$
d\log A={1\over A}dA - {1\over 2A{ }
^{2}}\left(\sum^{N}_{l=1}(\sum _{j=1}^N
c_{j}\sigma _{jl}\sum^{n}_{i=1}u_{ij})^{2}+\sum^{N}_{i=1}\beta _i^2\right)dt.
\eqno{(2.18)}
$$
By (2.1) and (2.4)
$$
\sum^{n}_{i=1}k_{i}dw_{i}=(\sum^{N}_{j=1}c_{j}
\mu _{j}\sum^{n}_{i=1}u_{ij})dt+\sum^{N}_{l=1}
(\sum^{N}_{j=1}c_{j}\sigma _{jl}\sum^{n}_{i=1}u_{ij})dB_{l}.
\eqno{(2.19)}
$$
Combining (2.18) and (2.19) one can get
$$
\sum^{n}_{i=1} A^{*}_{i}{dw{ } _{i}\over w{ } _{i}} +
d \log A - {d A\over A} ={1\over A}\sum^{N}_{l=1}(\sum^{N}_{j=1}
c_{j}\sigma _{jl}\sum^{n}_{i=1}u_{ij})dB_{l}+
\eqno{(2.20)}
$$
$$
+\left({1\over A}(\sum^{N}_{j=1}c_{j}\mu _{j}
\sum^{n}_{i=1}u_{ij}) - {1\over 2A{ } ^{2}}
\sum^{N}_{l=1}(\sum^{N}_{j=1}c_{j}\sigma _{jl}
\sum^{n}_{i=1}u_{ij})^{2}-\frac{1}{2A^2}\sum _{i=1}^N \beta _i^2\right) dt.
$$
In view of (2.16) and  It\^{o}'s formula,
$$
d\log k_i={dk_i\over k_i}-
{1\over 2k_i^2}\left(\frac{\beta _{i}}{w_i}\right)^2 dt.\eqno{(2.21)}
$$
Hence
$$
\sum _{i=1}^n {A^*_i}^2 {dk_i\over k_i}-
\sum _{i=1}^n {A^*_i}^2 d\log k_i=\frac{1}{2A^2}\sum _{i=1}^N\beta _i^2 dt.
\eqno{(2.22)}
$$
From (2.20) and (2.22) it follows
$$
\sum^{n}_{i=1} A^{*}_{i}{dw{ } _{i}
\over w{ } _{i}} +d\log A - {dA\over A}
+\sum _{i=1}^n {A^*_i}^2{dk_i\over k_i}-
\sum _{i=1}^n {A^*_i}^2 d\log k_i
$$
$$
=\sum^{N}_{l=1}{1\over A}(\sum^{N}_{j=1}
c_{j}\sigma _{jl}\sum^{n}_{i=1}u_{ij})dB_{l}+
$$
$$
+\left({1\over A}\sum^{N}_{j=1}c_{j}\mu _{j}
\sum^{n}_{i=1}u_{ij} - {1\over 2A{ } ^{2}}
\sum^{N}_{l=1}(\sum^{N}_{j=1}c_{j}\sigma _{jl}
\sum^{n}_{i=1}u_{ij})^{2}\right) dt.\eqno{(2.23)}
$$
By It\^{o}'s formula  for $f(x)=\exp(x)-1$ (see (2.8)) and (2.23) we have
$$
d\bar{r}(0,t)=(\bar{r}(0,t)+1)\frac{1}{A}\sum^{N}_{l=1}(\sum^{N}_{j=1}
c_{j}\sigma _{jl}\sum^{n}_{i=1}u_{ij})dB_{l}+
$$
$$+(\bar{r}(0,t)+1)\left({1\over A}\sum^{N}_{j=1}c
_{j}\mu _{j}\sum^{n}_{i=1}u_{ij} -
{1\over 2A{ } ^{2}}\sum^{N}_{l=1}(\sum^{N}_{j=1}
c_{j}\sigma _{jl}\sum^{n}_{i=1}u_{ij})^{2}+\right.
$$
$$
+\left. {1\over 2} \sum^{N}_{l=1}{1\over A{ } ^{2}}
(\sum^{N}_{j=1}c_{j}\sigma _{jl}\sum^{n}_{i=1}u_{ij})^{2}\right)dt.
$$
Hence
$$
d\bar{r}(0,t)=(\bar{r}(0,t)+1)
\frac{1}{A}\sum^{N}_{l=1}(\sum^{N}_{j=1}c_{j}\sigma _{jl}
\sum^{n}_{i=1}u_{ij})dB_{l}+
$$
$$
+(\bar{r}(0,t)+1){1\over A}
(\sum^{N}_{j=1}c_{j}\mu _{j}\sum^{n}_{i=1}u_{ij})dt.
$$
Since
$c_{i}(t)$ is an ${\Bbb F}$-martingale for each $i$,
we have $\mu _{i}(t)\equiv 0$ for each $i$,
and, in consequence,
$$
d(\bar{r}(0,t)+1)=(\bar{r}(0,t)+1)\frac{1}{A}
\sum^{N}_{l=1}(\sum^{N}_{j=1}c_{j}\sigma _{jl}\sum^{n}_{i=1}u_{ij})dB_{l}.
$$
By assumption of boundedness of $\sigma_{ij}$ we have
$${\Bbb E}\exp\left({1\over 2}\sum^{N}_{l=1}
 \int^{T}_{0}\Bigl(\frac{1}{A(t)}\sum^{N}_{j=1}\sum^{n}_{i=1}c_j(t)\sigma _{jl}(t)
u_{ij}(t)\Bigr)^{2}dt\right)<\infty,
$$
because $A(t)=\sum _{j=1}^N\sum _{i=1}^n c_j(t)u_{ij}(t)$.
Hence from the Novikov condition it follows that
$\{\bar{r}(0,t): 0\le t\le T\}$
is a martingale (see e.g.
Karatzas and Shreve(1998), p. 21,  or Shiryaev(1999)). $\Box$
\medskip

{\bf Remark. } If $c_i(\cdot)$ is a submartingale for each $i$, and
if short-selling is prohibited, then $\sum^{N}_{j=1}c_{j}\mu _{j}\sum^{n}_{i=1}u_{ij}\ge 0$ for every $t$ so
$\{\hat{r}(0,t): 0\le t\le T\}$ is a submartingale.

\medskip
Next, we formulate a list of properties which any properly defined
average rate of return should possess. Some of them can be found in
Kellison(1991). The average return
$\bar{r}$ meets the demands.

Put
$$ R(t)=\sum^{n}_{i=1}
\int^{t}_{s} A_i^*(u) {dw_{i}(u) \over w_{i}(u)} +\int^{t}_{s}d
\log A(u) - \int^{t}_{s}{dA(u)\over A(u)}+
$$
$$
+\sum _{i=1}^n\int _s^t {A^*_i}^2 (u){dk_i(u)\over k_i(u)}-
\sum _{i=1}^n\int _s^t {A^*_i} ^2(u)d\log k_i(u).
$$
\medskip

{\bf Property 1.} If the group consists only the $i$-th fund, then
$$
\bar{r}(s,t)={w_{i}(t)-w_{i}(s) \over w_{i}(s)}
$$

{\it Proof.}  Observe that
$$
dR=\frac{dw_i}{w_i}+d\log k_i +d \log w_i -\frac{d(w_ik_i)}{w_ik_i}+
$$
$$
+\frac{dk_i}{k_i}-d\log k_i=d\log w_i,\eqno{(2.24)}
$$
because of (2.15)-(2.16) and It\^{o}'s formula (2.9). From (2.24) we get
$$
\bar{r}(s,t)=\exp(\int _s^td\log w_i)-1=(w_i(t)-w_i(s))/w_i(s),
$$
and the proof is completed. $\Box$

\medskip
\medskip

{\bf Property 2.} For every $s<u<t$ with probability one,
$$
1+\bar{r}(s,t)=(1+\bar{r}(s,u))(1+\bar{r}(u,t))
$$

{\it Proof.} A direct consequence of the following property of It\^{o}'s integral: for
every $a<c<b$
$$
\int _a^b\xi dB=\int _a^c\xi dB+\int _c^b\xi dB, \quad {\Bbb P}-a.s. \Box
$$

\medskip

\medskip
{\bf Property 3.} If on a subset of a probability space the
accounting units of all funds have the same values
over $[s,t]$, i.e. $w_{1}(u)=w_{2}(u)=\ldots  =w_{n}(u)$
for every $u$, $s\le u\le t$, then the following equality
$$
\bar{r}(s,t)={w_{1}(t)-w_{1}(s)\over w_{1}(s)}.
$$
holds on the same subset.

{\it Proof.} Put $w =  w_1$ and put $k= \sum _{i=1}^n k_i$.
 Of course,
$$
A_i=w k_i, \quad A=w k.
$$
Hence
$$
dR=\sum _{i=1}^n A_i^*\frac{dw}{w} +d\log (wk)-
\frac{d(wk)}{wk}+\sum _{i=1}^n \frac{k_idk_i}{k^2}-\sum _{i=1}^n
\left(\frac{k_i}{k}\right)^2 d\log k_i
$$
$$
=d\log w + d\log k -\frac{dk}{k}+\sum _{i=1}^n \frac{k_idk_i}{k^2}-\sum _{i=1}^n
\left(\frac{k_i}{k}\right)^2 d\log k_i.\eqno{(2.25)}
$$
Since $w_1=\ldots =w_n=w$, we get from (2.2) and (2.7)
$$
dk=\frac{1}{w}\sum _{i=1}^n \alpha _i dt+
\sum _{i=1}^N \frac{\beta _i}{w}dB_{i+N}.
$$
By the local It\^{o} formula
$$
d\log k=\frac{dk}{k}-\frac{1}{2k^2}\sum _{i=1}^n\left(\frac{\beta _i}{w}\right)^2dt.
\eqno{(2.26)}
$$
From (2.21), we obtain
$$
\sum _{i=1}^n \frac{k_idk_i}{k^2}-\sum _{i=1}^n
\left(\frac{k_i}{k}\right)^2 d\log k_i=
\sum _{i=1}^n\left(\frac{k_i}{k}\right)^2\frac{1}{2k^2_i}\left(\frac{\beta _i}{w}\right)^2dt
$$
$$
=\frac{1}{2(kw)^2}\sum _{i=1}^n\beta _i^2 dt.\eqno{(2.27)}
$$
Combining (2.25)-(2.27) yields
$$
dR=d\log w,
$$
which completes the proof. $\Box $
\medskip

\medskip
{\bf Property 4.} Assume that there are reals $\alpha _i> 0$ and a function
$\phi :[s,t]\rightarrow [0,\infty )$  such that
$\sum _{i=1}^n \alpha _i=1$ and $k_i(u)=\alpha _i \phi (u)$ for
all $u\in [s,t]$, $i=1,\ldots ,n$. Then
$$
\bar{r}(s,t)={\sum _{i=1}^n \alpha _i r_i w_i(s)\over
\sum _{i=1}^n \alpha _i w_i(s)},
$$
where $r_i=(w_i(t)-w_i(s))/w_i(s)$. Moreover, if we assume that the number
of units of every fund is constant over the time interval $[s,t]$, i.e.
$\phi (u_1)=\phi (u_2)$ for $s\le u_{1}<u_{2}\le t$, then $$ \bar{r}(s,t)=
{A(t)-A(s)\over A(s)}. $$

{\it Proof.} By assumption,
$$
\sum_{i=1}^n{A^*_i}^2\frac{dw_i}{w_i}=
\frac{d(\sum _{i=1}^n\alpha _iw_i)}{\sum _{i=1}^n\alpha _iw_i}.
\eqno{(2.28)}
$$
We conclude from (2.15) that
$$
d\sum _{i=1}^n \alpha _i w_i=
\left(\sum _{i=1}^n{\alpha _i\over k_i}
\sum _{j=1}^N c_j\mu _j u_{ij}
\right)dt+
\sum _{l=1}^N
\left(\sum _{j=1}^Nc_j\sigma _{jl} \sum _{i=1}^n
{\alpha _i\over k_i} u_{ij}\right)dB_l.\eqno{(2.29)}
$$
By (2.29) and the It\^{o} formula
$$
d\log (\sum _{i=1}^n \alpha _i w_i)=
\frac{d(\sum_{i=1}^n\alpha _i w_i)}{\sum_{i=1}^n\alpha _i w_i}-
\frac{1}{2(\sum_{i=1}^n\alpha _i w_i)^2}
\sum _{l=1}^N
\left(\sum _{j=1}^Nc_j\sigma _{jl} \sum _{i=1}^n
\frac{1}{\phi } u_{ij}\right)^2dt
$$
$$
=\frac{d(\sum_{i=1}^n\alpha _i w_i)}{\sum_{i=1}^n\alpha _i w_i}-
\frac{1}{2A^2}
\sum _{l=1}^N
\left(\sum _{j=1}^Nc_j\sigma _{jl} \sum _{i=1}^n
u_{ij}\right)^2dt.\eqno{(2.30)}
$$
Combining (2.28) and (2.30) yields
$$
\sum_{i=1}^n{A^*_i}^2\frac{dw_i}{w_i}=d\log (\sum _{i=1}^n \alpha _i w_i)+
\frac{1}{2A^2}
\sum _{l=1}^N
\left(\sum _{j=1}^Nc_j\sigma _{jl} \sum _{i=1}^n
u_{ij}\right)^2dt.\eqno{(2.31)}
$$
From (2.18) and (2.31) we get
$$
\sum_{i=1}^n{A^*_i}^2\frac{dw_i}{w_i}=d\log (\sum _{i=1}^n \alpha _i w_i)+
\frac{1}{A}dA-d\log A-\frac{1}{2A^2}\sum _{i=1}^N \beta _i^2dt.\eqno{(2.32)}
$$
By (2.22) and (2.32), one can obtain
$$
dR=d\log (\sum _{i=1}^n\alpha _i w_i).
$$
Consequently,
$$
\bar{r}(s,t)=\exp \left(\int _s^td\log (\sum _{i=1}^n\alpha _i  w_i(u))du
\right)-1 $$ $$ = {\sum _{i=1}^n \alpha _i w_i(t)-\sum _{i=1}^n \alpha
_i w_i(s)\over \sum _{i=1}^n \alpha _i w_i(s)}={\sum _{i=1}^n \alpha _i r_i
w_i(s)\over \sum _{i=1}^n \alpha _i w_i(s)}. $$ Moreover, if $\phi (u)=c$
for each $u$ with $c>0$, then $$ \bar{r}(s,t)={\sum _{i=1}^n c\alpha _i
w_i(t)-\sum _{i=1}^n c\alpha _i w_i(s)\over \sum _{i=1}^n c\alpha _i
w_i(s)}=\frac{A(t)-A(s)}{A(s)}. $$
This completes the proof. $\Box $

\end{document}